\documentclass[11pt]{article}

\input{amssym.def}
\input{amssym}
\usepackage{eufrak}

\newtheorem{theorem}{Theorem}

\newtheorem{definition}[theorem]{Definition}

\newtheorem{lemma}[theorem]{Lemma}
\newtheorem{proposition}[theorem]{Proposition}
\newtheorem{remark}[theorem]{Remark}

\def\beq{\begin{equation}}
\def\be{\begin{equation}}
\def\eeq{\end{equation}}
\def\ee{\end{equation}}

\def\C{{\Bbb C}}
\def\R{{\Bbb R}}

\def\DD{\mathfrak{D}}
\def\Diff{{\mathfrak{Diff}}}
\def\VV{{{\cal{V}}}}
\def\Iso{{\rm{Iso}}}
\def\parh{\pa_{\rh}}

\def\D{{\frak D}}

\def\D{{\cal D}}

\def\pat{\partial_\tau}
\def\tpat{\tilde{\partial_\tau}}

\def\gg{\mathfrak{g}}
\def\pa{\partial}
\def\dd{\tilde{\partial}_t}

\def\tpa{\tilde{\pa_t}}
\def\V{\cal V}
\def\W{{\cal W}}
\def\K{\Bbb K}
\def\Z{\Bbb Z}

\def\ot{\otimes}
\def\h{h}

\def\hh{\hbar}
\def\rh{r_{\hh}}

\def\gh{\mathfrak{g}_h}

\def\Sym{{\rm Sym}}

\def\End{{\rm End}}

\def\De{{\Delta}}
\def\de{{\delta}}

\def\Cas{{\rm Cas}}
\def\Cs{{\rm Cas}}

\def\al{\alpha}

\def\Om{{\Omega}}

\def\la{{\lambda}}
\def\La{{\Lambda}}

\def\Ah{{\cal A}_\hh}

\def\Om{\Omega}

\begin{document}

\makeatletter
\renewcommand{\theequation}{{\thesection}.{\arabic{equation}}}
\@addtoreset{equation}{section}
\makeatother

\title{Noncommutative Geometry and dynamical models on $U(u(2))$ background}
\author{
\rule{0pt}{7mm} Dimitri
Gurevich\thanks{gurevich@univ-valenciennes.fr}\\
{\small\it LAMAV, Universit\'e de Valenciennes,
59313 Valenciennes, France}\\
\rule{0pt}{7mm} Pavel Saponov\thanks{Pavel.Saponov@ihep.ru}\\
{\small\it
National Research University Higher School of Economics,}\\
{\small\it
International Laboratory of Representation Theory and Mathematical Physics,}\\
{\small\it 20 Myasnitskaya Ulitsa, Moscow 101000, Russia}\\
{\footnotesize\it \&}\\
{\small\it Division of Theoretical Physics, IHEP, 142281
Protvino, Russia} }

\maketitle

\begin{abstract}
In our previous publications we have introduced  a differential calculus on the algebra $U(gl(m))$ based on a new
form of the Leibniz rule which differs from that usually employed in Noncommutative Geometry. This differential
calculus includes partial derivatives in generators of the algebra $U(gl(m))$ and their differentials. The corresponding
differential algebra $\Om(U(gl(m)))$ is a deformation of the commutative algebra $\Om(\Sym(gl(m)))$. A similar
claim is valid for the  Weyl algebra $\W(U(gl(m)))$ generated by the algebra $U(gl(m))$ and  the mentioned  partial
derivatives. In the particular case $m=2$ we treat the compact form $U(u(2))$ of this algebra as a quantization of the
Minkowski space algebra. Below, we consider noncommutative versions of the Klein-Gordon equation and the
Schr\"odinger equation for the hydrogen atom. To this end we define an extension  of the algebra $U(u(2))$
by adding to it meromorphic functions in the so-called quantum radius and quantum time.  For the quantum Klein-Gordon
model we get (under an assumption on momenta) an analog of the plane wave, for the quantum hydrogen atom
model we find the first order corrections to the ground state energy and the wave function.
\end{abstract}

{\bf AMS Mathematics Subject Classification, 2010:} 46L65, 46L87, 81T75


\section{Introduction}

One of the basic problems of the Noncommutative (NC) Geometry is to construct a differential calculus on a given NC algebra $A$. This
calculus should involve the differential algebra $\Om(A)$, generated by $A$ and  differentials $da, \, a\in A$, and vector fields (in appropriate cases,
analogs of partial derivatives). The core of the problem is a proper definition of the Leibniz rule for these vector fields and the corresponding
de Rham operator. Especially, this problem is interesting for algebras arising from quantization of commutative ones since in this case the
question of deformation properties of the  algebra  $\Om(A)$ can be  addressed.

The most common approach to this problem is based on the so-called K\"{a}hler forms which generate the  universal differential algebra. In this
algebra the usual Leibniz rule $d(ab)=da\, b+a\, db$ is preserved  but the the commutativity relation $a\, db=db\,a$ is canceled. As a result,
the differential algebra becomes bigger than the usual one, even if $A$ is a commutative algebra, say the coordinate algebra of a regular algebraic
variety\footnote{In order to reduce this differential algebra one sometimes passes to the quotient $\overline{\Om(A)}=\Om(A)/[\Om(A),\Om(A)]$
which is not an algebra at all. However, this passage leads to a change of the initial algebra $A$: it is replaced by the commutative quotient $A/[A,A]$.
Thus, the finial object $\overline{\Om(A)}$ cannot be treated as a differential algebra over the initial algebra $A$.}.

Sometimes, the Hochschild homology $H_i(A)$ is considered as a convenient analog of the space  $\Om^i(A)$ (for the relations between these two objects
for NC algebras $A$ we refer the reader to \cite{Lo}). Whereas the Hochschild
cohomology $H^i(A)$ plays the role of the space of polyvector fields on $A$ (this is motivated by the famous Hochschild-Kostant-Rosenberg theorem).
In particular, the first group $H^1(A)$ is identified with the space of derivations (quotiented over the inner ones), i.e. operators subject to the usual Leibniz
rule. An advanced calculus based on this treatment of the Hochschild (co)homology can be found in the papers \cite{TT, Gi, V} and others (also, see \cite{Lo, L}). Nevertheless, if $A$ is a coordinate algebra of a regular variety, the space $H_i(A)$ differs from that
$\Om^i(A)$ considered in its classical sense.

The same is valid if we replace $H_i(A)$ here by the cyclic homology $HC_i(A)$ due to
A.Connes. He also introduced another complex in which the
role of  the de Rham operator  is played by the commutator with a Dirac operator. However, it differs from  the usual de Rham complex as well.

Also, observe that in all mentioned approaches no analogs of partial derivatives
are introduced, even if a given NC algebra is a deformation of the coordinate algebra
of a linear space\footnote{The only deformation  algebras, on which partial 
derivatives can be introduced in the classical way are these arising from the Moyal quantization  of the Poisson brackets $\{x_i,x_j\}=\theta_{ij}$ with constant
$\theta_{ij}$.}.

A couple of years ago we suggested a new type of differential calculus (\cite{GPS2}, see also \cite{GS2}) on the algebras $U(gl(m))$, the super-algebras
$U(gl(m|n))$, and  their braided analogs --- the so-called modified Reflection Equation algebra ${\cal L}(R)$ related to a Hecke symmetry\footnote{The
Hecke symmetry is a solution of the quantum Yang-Baxter equation subject to a second degree polynomial relation.} $R$. The calculus is based on so-called {\em permutation relations} between elements of the NC algebra in question and
{\em partial derivatives} on it.

If such a Hecke symmetry $R$ is a deformation of the usual
flip $P$, the corresponding Reflection Equation algebra ${\cal L}(R)$ admits a construction of the differential algebra which is a deformation of its classical
counterpart. Moreover, the corresponding Weyl algebra\footnote{In physical literature such an algebra is usually called the Heisenberg algebra.}, generated by ${\cal L}(R)$ and the mentioned partial derivatives, is also a deformation of its classical
counterpart. Besides, if $R$ is involutive ($R^2=I$) the mentioned permutation relations can be expressed via a new form of the Leibniz rule which in its turn is related to a coproduct  in the
 algebra generated by the partial derivatives.

Below we  assume $R$ to be the usual flip $P$.
Moreover, we deal with a particular case $m=2$, i.e. our basic NC algebra
is\footnote{The parameter $\h$ in the notation $\gh$ means that we have introduced it as a factor at the  bracket of a Lie algebra $\gg$.
Note, that $h$ is not a formal parameter and it can be specialized.} $U(gl(2)_\h)$ or rather  its compact counterpart $U(u(2)_\h)$. The latter algebra can be treated as
a NC (quantum) analog of the Minkowski space algebra. Emphasize that the algebra $U(u(2)_\h)$ is a deformation (quantization) of the commutative algebra $\Sym(u(2))$
which is considered as the Minkowski space algebra provided it is endowed with a proper metric.

In this connection we  want to mention other deformations of the Minkowski space algebra such as q-Minkowski space algebra (in fact it is a particular case of the Reflection Equation algebra),
$\kappa$-Minkowski space algebra, and the algebras  mentioned in the footnote 2. There is a large number of publications
 devoted to these NC versions of the  Minkowski space algebra. Note, that our version of the quantum Minkowski space algebra, namely $U(u(2)_\h)$ has a lot of
advantages with respect to the other ones.

First of all, the algebra $U(u(2)_\h)$ is an $SO(3)$-covariant enveloping algebra of the reductive algebra $u(2)_\h$ whose simple component is the Lie algebra $su(2)_\h$. In virtue of the
Kostant's theorem (see section 2) the algebra $U(su(2)_\h)$ admits a decomposition in the direct sum of isotypic components
\be
\bigoplus_{k=0} \left( Z(U(su(2)_\h))\ot V^k\right)
\label{compo}
\ee
where $Z(U(su(2)_\h))$ is the center of the algebra  $U(su(2)_\h)$ and $V^k\subset U(su(2)_\h) $ is the irreducible $SU(2)$-module ($\dim V^k=2k+1$). This decomposition provides us with a basis in the algebra $U(su(2)_\h)$ convenient for our computations.

Besides, in the algebra $U(u(2)_\h)$ there is a version of the Cayley-Hamilton (CH) identity for a matrix $N$ composed of the $U(u(2)_\h)$ generators.
The roots $\mu_1$ and $\mu_2$ of this CH identity ("eigenvalues" of the  matrix $N$) are related to a quantity which can be treated as  a quantum analog of the
radius $\sqrt{x^2+y^2+z^2}$. This "quantum radius" $\rh=\frac{1}{2i}(\mu_1-\mu_2)$ together with the quantity $\tau=\frac{\mu_1+\mu_2-h}{2i}$ (treated as a quantum  analog
of the time) is very useful for a parametrization of the center $Z(U(su(2)_\h))$. It also enables us to construct a central extension of the algebra $U(su(2)_\h)$ (and
consequently, of $U(u(2)_\h)$) needed for constructing  NC versions of certain physical models.

In the present paper we are dealing with the model of the free massive scalar field described by the Klein-Gordon equation and
with the model of a particle in a central potential field (the hydrogen atom) described by the Schr\"odinger equation. In our previous publication we have considered the Klein-Gordon equation in nontrivial rotationally
invariant metric. Below we study the simplest form of the NC Klein-Gordon equation and find its plane wave type solution. For this solution we exhibit (under an assumption on momenta) a quantum
deformation of the so-called relativistic dispersion relation $E^2=m^2 c^4+p^2c^2$ connecting the energy $E$ of a particle of the mass $m$ and the square of its
momentum $p^2$.  As for the NC Schr\"odinger equation, we find the first order corrections to the ground state energy and to the wave function of the corresponding
hydrogen atom  model.

Emphasize that there are known numerous attempts to generalize some models of mathematical physics to quantum algebras. Besides  the aforementioned versions of
the Minkowski space algebra, we want to cite the papers \cite{Ku1, Ku2} (also, see the references therein) where some rotationally invariant models were quantized.
The essential  dissimilarity between the cited papers and ours consists in different treatment of the differential operators (such as the Laplacian $\Delta = \pa_x^2+\pa_y^2+\pa_y^2$)
appearing in dynamical equations of the models\footnote{Some papers deal with "fuzzy varieties" where the role of the Laplacian is played by the image of the Casimir element,
see \cite{CMS}.}. In the cited papers the partial derivatives are treated in the classical sense, only the usual
commutative product of  functions is replaced by its deformed version (called $\star$-product). By contrary, in our
approach the partial derivatives are directly defined on the NC algebra in question
and their  properties differ drastically from those of the usual derivatives.

First, the partial derivatives in the spacial generators $x,y,z$ (and consequently,  the Laplacian) act
 nontrivially  on the time variable and visa versa.

Second, we treat the quantized models in terms of the quantum algebra without using any $\star$-product. The only problem arising in our approach is  pulling an operator
which describes the model forward to the quantum algebra. We call this procedure {\em the quantization of operators}. For the operators we are
dealing with, this quantization is easy: the derivatives should be treated in the new sense (i.e. according to our definition of the derivatives  on the algebra $U(u(2)_h)$) and
the radial variable $r$ in a potential $V(r)$  should be replaced by its quantum counterpart $r_\hbar$, as well as the time variable $t$ is replaced by its quantum counterpart
$\tau$.  In more detail the problem of quantization of differential operators is considered in \cite{GS2} (also,
see section 6 below.)

Third, the models in question can be restricted to any component coming  in the sum (\ref{compo}) and then to some lattices in 2-dimensional space  "time-quantum radius". This property is due to the fact that the partial derivatives in time and quantum radius are difference
operators\footnote{Note, that these operators are "additive" compared with "multiplicative" $q$-derivative appearing in connection with Quantum Group theory.}.
In this sense we speak about a discretization of the space-time.

The ground field $\K$ is assumed to be $\R$ or $\C$,  its concrete form is always clear  from the context.

The paper is organized as follows. In the next section we  introduce a quantum version of the Weyl algebra. Also, we define its extended version by adding meromorphic
functions in the quantum radius and the time to it. In the section 3 we give an operator meaning to elements of the quantum Weyl algebra --- they become operators
acting on the algebra $U(u(2)_\h)$. In section 4 we introduce the partial derivative in quantum radius and clarify its role in discretization of models we are dealing with.
Section 5 is devoted to our quantization of the aforementioned rotationally invariant dynamical models. In section 6 we present some concluding remarks and discuss the
perspective of our study. In Appendix we get together some technical results used in the text.

\vspace*{3 mm}

\noindent
{\bf Acknowledgement}

This study (research grant No 14-01-0173) was supported by The National Research University Higher School of Economics Academic Fund Program in 2014/2015.

\section{The Weyl algebra $\W(U(u(2)_\h))$ and its extension}

First, recall that according to \cite{GPS2, GS2} the Weyl algebra $\W(U(gl(m)_\h))$ is a unital associative algebra generated by
$2 m^2$ elements $\{n_i^j\}_{1\le i,j\le m}$ and $\{d_i^j\}_{1\le i,j\le m}$. Having defined the $m\times m$ {\it generating
matrices} $N=\|n_i^j\|$ and $D=\|d_i^j\|$, we present the multiplication rules of the  $\W(U(gl(m)_\h))$ generators in the
following matrix form:
\be
\begin{array}{rcl}
PN_1PN_1-N_1PN_1P&=&h (PN_1- N_1P), \\
PD_1PD_1&=&D_1PD_1P,\\
D_1PN_1 P-PN_1P\,D_1&=&P+\h\,D_1P,
\end{array}
\label{sett}
\ee
where $h$ is a numeric parameter, $N_1 = N\otimes I_{m\times m}$, $D_1  = D\otimes I_{m\times m}$ and $m^2\times m^2$
matrix $P$ is the usual flip.

As is clear from the first line of the system (\ref{sett}), $n_i^j$ are generators of the subalgebra $U(gl(m)_\h)\subset \W(U(gl(m)_\h))$.
The second line expresses the commutativity of generators $d_i^j$. The corresponding Abelian subalgebra will be denoted
$\cal D$. The relations coming in the third line are called {\em the permutation relations}. The permutation relations allow us
to treat the generators $\{d_i^j\}$ as (quantum) partial derivatives. In particular, in the classical limit ($\h=0$) the element
$d_i^j$ can be identified with the partial derivative $\partial/\partial{n_j^i}$. From now on we shall use the term ``partial derivatives'' for
the elements $d_i^j$, while the precise operator meaning will be given to them in the next section.

Note that the algebra $\W(U(gl(m)_\h))$ can be treated as the enveloping algebra of a semidirect product of the Lie algebra $gl(m)_\h$
with the basis $n_i^j$ and the commutative algebra ${\cal D} = gl(m)_0$ with the basis $d_i^j$. In fact we have a slight modification of the
semidirect product because of the constant terms, entering the flip $P$ in the third line of (\ref{sett}), but this does not impact the claim that in virtue of the PBW theorem the family of elements
$$
(n_1^1)^{a_1^1}(n_1^2)^{a_1^2}\dots (n_m^m)^{a_m^m}(d_1^1)^{b_1^1}(d_1^2)^{b_1^2}\dots (d_m^m)^{b_m^m},\qquad a_i^j,\, b_i^j\in \Bbb{N}
$$
is a basis of the algebra $\W(U(gl(m)_\h))$.

\begin{remark}\rm  S.Meljanac and Z.{\v{S}}koda have suggested the realization of relations  (\ref{sett}) analogous to their construction
from \cite{MS}  (also, see \cite{S}). Let $\tilde{N} = \|t_i^j\|$ be the generating matrix of the algebra $\Sym(gl(m))$ (i.e. it satisfies the first
line relations (\ref{sett}) with $h=0$) and $D = \|d_i^j\|$ be the matrix composed of the usual partial derivatives in $t_i^j$: $d_i^j =
\partial/\partial t_j^i$. Then the matrices $N=\tilde{N}+h \tilde{N} D$ and $D$ satisfy the system (\ref{sett}).
\end{remark}

Now, we write down the relations (\ref{sett}) in the particular case $m=2$ and pass to
the compact counterpart $U(u(2)_\h)$ of the algebra $U(gl(2)_\h)$.
To simplify formulae we use the following notation for the $U(gl(2)_\h)$ generators:
$$
n_1^1=a,\qquad n_1^2=b,\qquad n_2^1=c,\qquad n_2^2=d.
$$
The generators  of the algebra $U(u(2)_\h)$
\be t=\frac{a+d}{2},\qquad x=\frac{i(b+c)}{2},\qquad y=\frac{c-b}{2},\qquad z=\frac{i(a-d)}{2}.
\label{abcd}
\ee
are subject to the system of relations
$$
[x, \, y]=\h\, z,\quad [y, \, z]=\h\, x,\quad[z, \, x]=\h\, y,\quad [t, \, x]=[t, \, y]=[t, \, z]=0.
$$
The center $Z(U(u(2)_\h))$ is generated by the element $t$ and quadratic Casimir element
$$
\Cas = x^2+y^2+z^2.
$$

The partial derivatives in $t$, $x$, $y$ and $z$ will be denoted as $\partial_t$, $\partial_x$, $\partial_y$ and
$\partial_z$ respectively. Besides, it is convenient to introduce a shifted derivative $\tilde \partial_t = \partial_t +
\frac{2}{h}I$.

In terms of generators $t$, $x$, $y$, $z$ and the corresponding partial derivatives the third line of the system (\ref{sett}) reads
\be
\begin{array}{l@{\quad}l@{\quad}l@{\quad}l}
\tilde\pa_t\,t - t\,\tilde\pa_t = {h\over 2}\,\tilde\pa_t & \tilde\pa_t\, x - x\,\tilde\pa_t
=-{h\over 2}\,\pa_x &
\tilde\pa_t\, y - y\, \tilde\pa_t=-{h\over 2}\,\pa_y &\tilde\pa_t\, z - z\,\tilde\pa_t=- {h\over 2}\,\pa_z\\
\rule{0pt}{7mm}
\pa_x\, t - t\,\pa_x = {h\over 2}\,\pa_x &\pa_x \,x -  x\,\pa_x = {h\over 2}\,\tilde\pa_t &
\pa_x \, y-  y\,\pa_x = {h\over 2}\,\pa_z & \pa_x \,z - z\, \pa_x  = - {h\over 2}\,\pa_y \\
\rule{0pt}{7mm}
\pa_y \,t - t \, \pa_y = {h\over 2}\,\pa_y & \pa_y \,x -  x\,  \pa_y = -{h\over 2}\,\pa_z &
\pa_y \,y - y \,  \pa_y = {h\over 2}\,\tilde\pa_t & \pa_y \,z - z \,  \pa_ y= {h\over 2}\,\pa_x\\
\rule{0pt}{7mm}
\pa_z \,t - t \,\pa_z = {h\over 2}\,\pa_z & \pa_z \,x - x \,\pa_z = {h\over 2}\,\pa_y&
\pa_z \,y -  y\,\pa_z = -{h\over 2}\,\pa_x & \pa_z \,z - z \,\pa_z = {h\over 2}\,\tilde\pa_t.
\end{array}
\label{leib-r}
\ee

Our next aim is to introduce quantum analogs of "radial variables". To this end we need a version of the
Cayley-Hamilton (CH) identity for the generating matrix $N$
\be
N=\left(\begin{array}{cc}
a& b\\
c& d
\end{array}\right)=\left(\begin{array}{cc}
t-i\,z& -i\, x-y\\
-i\, x+y& t+i\, z
\end{array}\right).
\label{matrN}
\ee

In this case  the corresponding characteristic polynomial
$\chi(\lambda)$ is quadratic and the CH identity explicitly reads
\be
\chi(N) = N^2-(2\,t+\h)\, N+\,(t^2+x^2+y^2+z^2+\h\, t)\,I= 0.
\label{CH}
\ee
Note, that the coefficients of the above identity are central in the algebra $U(u(2)_\h)$. Recall, that this is also true in the general case $U(gl(m)_\h)$.

The roots $\mu_1$ and $\mu_2$ of the characteristic polynomial are called the eigenvalues of the matrix $N$. This means that they solve the
system of two algebraic equations
\be
\mu_1+\mu_2=2\,t+\h ,\qquad
\mu_1 \,\mu_2= t^2+x^2+y^2+z^2+\h\, t,
\label{mmu}
\ee
and, consequently, the eigenvalues should be treated as elements of an algebraic extension of the center  $Z(U(u(2)_\h))$.

Equations (\ref{mmu}) allow us to express the center generators (and, therefore, any element of $Z(U(u(2)_\h))$) in terms of $\mu_1$ and $\mu_2$:
\be
t=\frac{\mu_1+\mu_2-\h}{2},\qquad \Cas=\frac{\h^2-\mu^2}{4}.
\label{cas-mm}
\ee
Hereafter, we use the notation $\mu=\mu_1-\mu_2$.

\begin{remark} \label{rem2} \rm Note that the quantity $\mu$ is defined up to a sign since the ordering of the eigenvalues is arbitrary.
\end{remark}

Now, we compute the permutation relations of the partial derivatives $\pa_x$, $\pa_y$, $\pa_z$ and the Casimir element.
In contrast with the commutative case, the other elements of the algebra $\W(U(u(2)_\h))$ turn out to be involved in these permutation
relations. We need the following operators
\be
Q=x\pa_x+y\pa_y+z\pa_z,\quad X=y\pa_z-z\pa_y,\quad Y=z\pa_x-x\pa_z,\quad Z=x\pa_y-y\pa_x.
\label{QX}
\ee

\begin{theorem}
\be
\left(\!\!\!\begin{array}{c}
\partial_x\\
\rule{0pt}{6mm}
x\,\tilde\partial_t\\
\rule{0pt}{6mm}
x\,Q\\
\rule{0pt}{6mm}
X
\end{array}\!\!\!\right) \Cs =
\left(\!\!
\begin{array}{cccc}
\displaystyle \Cs -\frac{3}{4}\,\h^2&\h&0&\h\\
0&\displaystyle \Cs -\frac{3}{4}\,\h^2&-\h&0\\
0 &\h\,\Cs&\displaystyle \Cs + \frac{\h^2}{4} & 0\\
-\h\,\Cs &\h^2& \h&\displaystyle \Cs + \frac{\h^2}{4}
\end{array}
\!\!\right)
\left(\!\!\!\begin{array}{c}
\partial_x\\
\rule{0pt}{6mm}
x\,\tilde\partial_t\\
\rule{0pt}{6mm}
x\,Q\\
\rule{0pt}{6mm}
X
\end{array}\!\!\!\right) .
\label{form}
\ee
The permutation relation of $\partial_y$ (respectively, $\partial_z$) with the Casimir element are obtained from  (\ref{form}) by replacing
$x$ and $X$ for $y$ and $Y$ (respectively, for $z$ and $Z$).
\end{theorem}

We denote the $4\times 4$ matrix in  formula (\ref{form}) by $\Psi(\Cs)$. Being expressed in terms of $\mu$ (see (\ref{cas-mm})) it reads
\be
\Psi(\Cs)=\left(\!\!
\begin{array}{cccc}
\displaystyle \frac{-\mu^2-2\h^2}{4}&\h&0&\h\\
0&\displaystyle \frac{-\mu^2-2\h^2}{4}&-\h&0\\
0 &\h\,\frac{\h^2-\mu^2}{4}&\displaystyle \frac{2\h^2-\mu^2}{4} & 0\\
-\h\,\frac{\h^2-\mu^2}{4} &\h^2& \h&\displaystyle \frac{2\h^2-\mu^2}{4}
\end{array}
\!\!\right)
\ee

The matrix $\Psi(\Cs)$ is semisimple with the following spectrum:
$$
\Psi(\Cs)\sim \mathrm{diag}(\lambda_1,\lambda_1,\lambda_2,\lambda_2),\qquad
\lambda_1 = -\,\frac{\mu^2-2\h\mu}{4}, \quad
\lambda_2 = -\,\frac{\mu^2+2\h\mu}{4} .
$$
Its spectral decomposition reads
$$
\Psi(\Cs) =\lambda_1P_1(\mu) +\lambda_2P_2(\mu),\qquad P_1(\mu) =
\frac{\Psi(\Cs) - \lambda_2 I}{\lambda_1-\lambda_2},\quad
P_2(\mu) = \frac{\Psi(\Cs) - \lambda_1 I}{\lambda_2-\lambda_1},
$$
where $I$ is the unit $4\times 4$ matrix. The explicit matrix form of the complementary projectors $P_i(\mu)$ is as follows:
$$
P_1(\mu) =
\frac{1}{\mu}\,\left(\!\!
\begin{array}{cccc}
 \frac{\mu-\h}{2}&1&0&1\\
\rule{0pt}{6mm}
0&  \frac{\mu-\h}{2}&-1&0\\
\rule{0pt}{6mm}
0 & \frac{\h^2-\mu^2}{4} &  \frac{\mu+\h}{2}& 0\\
\rule{0pt}{6mm}
\frac{\mu^2-\h^2}{4}& \h& 1& \frac{\mu+\h}{2}
\end{array}
\!\!\right),\qquad
P_2(\mu) =
\frac{1}{\mu}\,\left(\!\!
\begin{array}{cccc}
 \frac{\mu+\h}{2}&-1&0&-1\\
\rule{0pt}{6mm}
0&  \frac{\mu+\h}{2}&1&0\\
\rule{0pt}{6mm}
0 & \frac{\mu^2-\h^2}{4} &  \frac{\mu-\h}{2}& 0\\
\rule{0pt}{6mm}
 \frac{\h^2-\mu^2}{4}&-\h& -1& \frac{\mu-\h}{2}
\end{array}
\!\!\right).
$$

The permutation relations of the column entering (\ref{form}) with the element $\Cas^p$, $p\ge 2$, can be obtained in the same way but
the matrix $\Psi(\Cas)$ must be replaced by its p-th degree. Using the spectral decomposition, we calculate $\Psi^p(\Cas)$ as
$$
\Psi^p(\Cas)=\la_1^p \,P_1(\mu)+\la_2^p \,P_2(\mu).
$$
So, we can find the permutation relations of the set of elements $\partial_x$, $x\tilde\partial_t$, $xQ$ and $X$ with  any polynomial  or even a series
$f(u)$ in $u=\Cas$. The only thing we need for this purpose is
$$
f(\Psi(\Cs))=f(\la_1) \,P_1(\mu)+f(\la_2) \,P_2(\mu).
$$

However, we want to get  the permutation relations of the above operators with $\mu$ and its powers $\mu^p, p\in \Z$.
Since $\mu = f(\Cas)$
where $f(u) = \sqrt{h^2-4u}$, we get
$$
f(\la_1)=\sqrt{(\mu-\h)^2}=\mu-\h,\,\,  f(\la_2)=\sqrt{(\mu+\h)^2}=\mu+\h.
$$
The sign is motivated by the classical limit. Finally, we arrive to the following permutation relations
\be
\left(\!\!\!\begin{array}{c}
\partial_x\\
\rule{0pt}{6mm}
x\,\tilde\partial_t\\
\rule{0pt}{6mm}
x\,Q\\
\rule{0pt}{6mm}
X
\end{array}\!\!\!\right)\mu^p=
(\mu-\h)^p\,P_1(\mu)
\left(\!\!\!\begin{array}{c}
\partial_x\\
\rule{0pt}{6mm}
x\,\tilde\partial_t\\
\rule{0pt}{6mm}
x\,Q\\
\rule{0pt}{6mm}
X
\end{array}\!\!\!\right)
+
(\mu+\h)^p\,P_2(\mu)
\left(\!\!\!\begin{array}{c}
\partial_x\\
\rule{0pt}{6mm}
x\,\tilde\partial_t\\
\rule{0pt}{6mm}
x\,Q\\
\rule{0pt}{6mm}
X
\end{array}\!\!\!\right).
\label{formm}
\ee

In the same way the permutation relations of the partial derivatives $\pa_y$ and $\pa_z$ with $\mu^p$, $ p\in \Z$ can be found.
As for the permutation relations of the operators $\tpa$ and $Q$ with $\mu^p$, they can be extracted from formula (\ref{formm}):
it suffices to cancel $x$ in the second and third lines of this matrix. Thus, we arrive to the result (also, see \cite{GS2} where these
formulae were first obtained)
\be
\begin{array}{l}
\displaystyle
 \tilde\partial_t \mu^p = \frac{(\mu+h)^p}{\mu}\left( Q+\frac{\mu+h}{2}\,\tilde\partial_t\right)-\frac{(\mu-h)^p}{\mu}\left(Q-\frac{\mu-h}{2}\,\tilde\partial_t \right)\\
\rule{0pt}{10mm}\displaystyle
Q\mu^p = \frac{(\mu+h)^p}{\mu}\frac{(\mu-h)}{2}\left( Q+\frac{\mu+h}{2}\,\tilde\partial_t\right)+\frac{(\mu-h)^p}{\mu}\frac{(\mu+h)}{2}\left(Q-\frac{\mu-h}{2}\,\tilde\partial_t \right).
\end{array}
\label{formmm}
\ee

Note that in the classical case ($h=0$) the variable $\mu$ is a multiple of the radius $r=\sqrt{x^2+y^2+z^2}$, namely,  $\mu=2i r$. In the general case
($h\not=0$) we introduce an analog of the radius by a similar formula.
\begin{definition} \label{def:4}\rm
The quantity $\displaystyle \rh=\frac{\mu}{2i}$ is called {\it the quantum radius}.
\end{definition}

Below, we use the element $\rh$ as well as $\mu$. However, dealing with $\rh$ we prefer to use the renormalized deformation parameter $\hbar = \frac{h}{2i}$
in order to avoid complex shifts of the type $\rh\to \rh+ia$ (as for the factor 2 it is introduced for the future convenience). For the same reason
we will use the renormalized  time $\tau=-it$. All our formulae represented in terms  of the initial quantities $(\mu, t, h)$ can be easily rewritten in terms of $(\rh, \tau, \hh)$
and visa versa. For example, the Casimir element in terms of new quantities has the form
$$
\Cas=\rh^2-\hh^2.
$$

We also define the partial derivative $\pat=i \pa_t$ and the corresponding "shifted derivative" $\tpat=i\tpa=\pat+\frac{1}{\hh}I$.
Note, that $[\partial_\tau, \tau] = \hbar\partial_\tau$.

Introduce now the commutative algebra $\K(\tau, r_\hbar) $ of functions in $\tau$ and $\rh$. For any  $f(\tau,\rh)\in\K(\tau,r_\hbar)$ we assume the function
$f(z_1,z_2)$ to be meromorphic functions in two complex variables. Next, consider the algebra
\be
\Ah=\left(U(su(2)_\h)\ot \K(\tau,\rh)\right)/\langle x^2+y^2+z^2-\rh^2+\hh^2 \rangle. \label{quot} \ee

The generators $\tau$ and $\rh$ are assumed to commute with any element of the algebra $\Ah$. Thus, the subalgebra $\K(\tau,r_\hbar)$ forms evidently the center $Z(\Ah)$ of this algebra.

Our next aim is to find the permutation relations of the partial derivatives $\tilde\partial_\tau,\dots,\pa_z$ and elements of the algebra $\K(\tau,\rh)$.

First, consider in $\K(\tau,r_\hbar)$ polynomial functions in two variables. Since $\tau$ and $r_\hbar$ commute, any polynomial $f(\tau,r_\hbar)$ can be presented
as a finite sum $f(\tau,r_\hbar)=\sum_i g_i(\tau)q_i(r_\hbar)$ where $g_i$ and $q_i$ are some polynomials in one variable. Then, as was shown
in \cite{GS2}, the permutation relations of the partial derivatives with a polynomial $g(\tau)$ are
\be
\tilde\partial_\tau g(\tau)=g(\tau+\hh)\tilde\partial_\tau,\qquad\partial_a g(\tau)=g(\tau+\hh)\partial_a,\quad a\in\{x,y,z\}.
\label{fo}
\ee
As for the  permutation relations  of the partial derivatives with $r_\hbar^p$ (and, consequently,
with any polynomial $f(r_\hbar)$), they  can be easily deduced  from  (\ref{formm}).

Therefore, the above considerations together with (\ref{fo}) and (\ref{formm}) allow us to get for any {\it polynomial} $f(\tau,r_\hbar)$ in two
variables the following result
\begin{eqnarray}
\pa_x f(\tau, \rh) =\frac{f(\tau+\hh, \rh+\hh)}{2\rh}\left((\rh \right.&\!\!\!\!\!+&\!\!\!\!\!\left.\hh)\pa_x+x\tilde\partial_\tau+iX\right)\nonumber\\
&\!\!\!\!\!+&\!\!\!\! \frac{f(\tau+\hh, \rh-\hh)}{2\rh}\left((\rh-\hh)\pa_x-x\tilde\partial_\tau-iX\right),
\label{dxf}
\end{eqnarray}
and similarly for the derivatives $\pa_y$ and $\pa_z$.

For the partial derivative $\tilde\partial_\tau$ we find
\be
\tilde\partial_\tau f(\tau, r_\hbar) =\frac{f(\tau+\hh, \rh+\hh)}{2\rh}\left((\rh+\hh)\tilde\partial_\tau+Q\right)+
\frac{f(\tau+\hh, \rh-\hh)}{2\rh}\left((\rh-\hh)\tilde\partial_\tau-Q\right).
\label{dtf}
\ee

We extend these formulae from polynomials to arbitrary functions $f(\tau,\rh)\in \K(\tau,\rh)$
by definition.

Now, define the  algebra $\W(\Ah)$ to be the algebra generated by  $\Ah$ and  the partial derivatives $\tilde\partial_\tau,\dots,\pa_z$ endowed with the
permutation relations described by formulae (\ref{leib-r}),  (\ref{dxf}) with its $y$- and $z$-analogs and (\ref{dtf}).
We call $\W(\Ah)$ {\em the extended Weyl algebra}.

Below, we define the derivative $\parh$ in $\W(\Ah)$  but it will be expressed in terms of already existing elements
of $\W(\Ah)$ so we shall not need its new extension.

By resuming this section, we want to emphasize  that the  algebra $\W(\Ah)$ is introduced as a quotient-algebra.
Intuitively, it is clear that it is a deformation of its classical counterpart, defined by formula (\ref{quot})
with $\h=0$. However, the fact that the algebra $\W(U(u(2)_h))$ is a deformation of its classical counterpart
$\W(\Sym(u(2)))$ can be shown rigorously. Indeed, it follows from the discussed above treatment of the algebra
$\W(U(gl(m)_h))$ as the enveloping algebra of a semidirect product.

\section{Differential  operators on the  algebra $\Ah$}

Now, we are going to define an action of the algebra  $\W(\Ah)$ onto $\Ah$. To this end we first define this
action for the  ``partial derivatives'' $d_i^j$ (recall that at the limit $\h=0$ we have $d_i^j=\pa_{n_j^i}$)
$$
{\cal D}\ot U(gl(m)_\h) \to U(gl(m)_\h):\quad  d_i^j\otimes n\mapsto \tilde n=d_i^j(n).
$$
Hereafter, the action of an operator $\DD$ on an element $n\in U(u(m)_h)$ will be denoted  $\DD(n)$, while we keep the notation
${\DD}n = \tilde n\tilde{\DD}$ for the permutation relation.

Below we need the counit map $\varepsilon:\D\to \K$ introduced by
\be
\varepsilon(1_\D)=1,\qquad \varepsilon(d_i^j)=0,\qquad \varepsilon(d_1d_2)=\varepsilon(d_1)\varepsilon(d_2),\quad
 \forall\,\,d_1,d_2\in \D.
\label{coun}
\ee
The symbol $1_{\cal A}$ stands for the unit element of the (sub)algebra $\cal A$.

Then we proceed as follows. Given a product of elements $d\,n$ (hereafter, we  omit the sign $\ot$), we first apply
the permutation relations (the third line of (\ref{sett})) to the elements $d$ and $n$ in order to move elements of $\D$
to the most right position
$$
d\,n = \sum_i n'_i d'_i +n'' 1_{\cal D},\qquad  n'_i, n''\in U(gl(m)_h),\quad d_i'\in {\cal D},
$$
where in the right hand side we  have  a finite sum of elements.
After that, we apply the counit map $\varepsilon$ to all elements of ${\cal D}$ in the above sum and get an element from $U(gl(m)_h)$:
$$
\sum_i n'_i d'_i + n''1_{\cal D}\stackrel{{\rm id\otimes\varepsilon}}{\longrightarrow} \sum_i n'_i\varepsilon(d'_i) +n''\varepsilon(1_{\cal D}) = n'',
$$
where we identify $U(gl(m)_h)\ot \K$ with $U(gl(m)_h)$.

Thus, all partial derivatives $\pa_t,...,\pa_z$ become operators acting on $U(u(2)_h)$.
Emphasize that the action of the partial derivatives is not subject to the classical Leibniz rule.
By contrary, these derivatives meet a new version of the Leibniz rule which was found by S.Meljanac and Z. \v{S}koda
(see also \cite{GS2} where we have given another but equivalent form of this Leibniz rule).
In order to describe this form of the Leibniz rule, consider the following coproduct in the algebra ${\cal D}$
\be
\De(d_i^j)=d_i^j\ot 1+1\ot d_i^j+h\sum_k d^j_k\ot d^k_i. \label{cco}
\ee

This coproduct  and  the counit (\ref{coun}) endow ${\cal D}$ with a  structure of a coassociative  coalgebra. Moreover, since
this structure is coordinated with the algebraic structure of ${\cal D}$, this algebra becomes  a bi-algebra.

Now, the standard action $d_i^j(n_k^l)=\de^l_i \de_k^j$ of the derivatives $d_i^j$ on the generators $n_k^l$  of the algebra
$U(gl(m)_h)$ can be extended to the whole algebra $U(gl(m)_h)$ via the above coproduct. In order to prove that the extended
action respects the algebraic structure of $U(gl(m)_h)$, it suffices to check that the ideal generated in this algebra by the elements
$$
P N_1PN_1-N_1 P N_1 P -h(P N_1-N_1 P)
$$
is invariant with respect to the action of the derivatives. The detail is left to the reader.

\begin{remark}{\rm Note that by using the method of \cite{GPS1} it is possible to recover the permutation relations (the third line of (\ref{sett})) on the base of the coproduct (\ref{cco}). However, in a more general case related to the Reflection Equation algebra (see
\cite{GS2}) we have not succeeded in finding the coproduct
corresponding to the permutation relations studied in this paper. Thus,  the permutation relations are in a sense a more fundamental structure enabling us to
define an action of the algebra $\cal D$ on that $U(gl(m)_h)$ without any coproduct.}
\end{remark}

Now, pass to the case $m=2$ and  extend the derivatives $\pa_t,...,\pa_z$ onto the  algebra $\Ah$
(here we deal with the compact form of the algebra $gl(2)_h$ and also we use the generator $\tau$ instead of $t$).
First, we exhibit the coproduct $\De$ in the generators $\pa_t, \pa_x,\pa_y, \pa_z$. Thus, we have
$$
\De(\pa_t)=\pa_t\ot 1 +1\ot \pa_t+\frac{h}{2}(\pa_t\ot \pa_t-\pa_x\ot \pa_x-\pa_y\ot \pa_y-\pa_z\ot \pa_z),
$$
$$
\De(\pa_x)=\pa_x\ot 1 +1\ot \pa_x+\frac{h}{2}(\pa_t\ot \pa_x+\pa_x\ot \pa_t+\pa_y\ot \pa_z-\pa_z\ot \pa_y).
$$
By applying the cyclic permutations to the latter formula, we get  $\De(\pa_y)$ and $\De(\pa_z)$.

Now, using formulae (\ref{dtf}) and (\ref{dxf}) we can compute the action of the derivatives
$\partial_\tau=i\pa_t$ and $\partial_x$ on
$f(\tau, r_\hbar)$:
\be
\pat(f(\tau, \rh))=\frac{f(\tau+\hh,\rh+\hh)(\rh+\hh)+f(\tau+\hh,\rh-\hh)(\rh-\hh)-2\rh f(\tau, \rh)}{2\rh\hh},
\label{patt}
\ee
\be
\pa_x(f(\tau, \rh))=\frac{x}{\rh}\frac{(f(\tau+\hh,\rh+\hh)-f(\tau+\hh,\rh-\hh))}{2\hh}.
\label{pax}
\ee
Similar formulae  are valid for  $\pa_y(f(\tau, \rh))$ or $\pa_z(f(\tau, \rh))$.

Extend the action of the derivatives on the algebra $\Ah$ via the above coproduct. In order to show that
the extended derivatives  are well defined on this algebra it sufficed to check that they map the ideal
generated in this algebra  by the elements
\be u g(\tau, \rh)-g(\tau, \rh)u,\,\, u\in\{\tau, x, y, z\} \label{elem}  \ee
to itself. The following proposition entails this property.

\begin{proposition} Each element
\be \pa_v(u f(\tau, \rh)-f(\tau, \rh)u),\,\, u,v \in\{\tau, x, y, z\}  \label{that} \ee
is a linear combinations of the elements of the form (\ref{elem}).
\end{proposition}

{\bf Proof.}  In fact, it suffices to check  (\ref{that}) with $v=x$ and $u=y$. Other relations are trivial or can be obtained by the
cyclic permutation of the generators $x,y, z$. Checking (\ref{that}) with $v=x,\,u=y$
can be done with the help of formula (\ref{pax}) or its equivalent form (\ref{xcase}) below.
\hfill\rule{6.5pt}{6.5pt}

Note that any element of the algebra $\W(U(u(2)_\h))$ can be presented as a finite sum of elements
$\pa_z^d\,\pa_y^c\,\pa_z^b\, \pa_t^a$, $a,b,c,d\in \Bbb{N}$
with coefficients from $U(u(2)_\h)$ (we put them on the left hand side). In a similar way we can realize
elements of the algebra  $\W(\Ah)$. This realization of any element of the algebra  ${\W}(U(u(2)_\h))$ or $\W(\Ah)$ is called
{\em canonical}. A canonical form of a given element $A\in \W(\Ah)$ allows us to assign to $A$ an operator
$\DD_A: \Ah\rightarrow \Ah$ (for this assignment we use the operator meaning of elements from ${\cal D}$ discussed above).

Let us denote $\Diff$ the vector space of all such differential operators. Thus, we have an isomorphism
$$
\Iso:\W(\Ah)\to \Diff,
$$
which is a representation of the algebra $\W(\Ah)$.

Now we construct a special basis in the algebra $U(u(2)_\h)$ useful in the sequel.
First, note that $U(u(2)_\h)=U(su(2)_\h)\ot \R[t]$. Second, we use the fact that the algebra $U(su(2)_\h)$ is free
over its center $Z(U(su(2)_\h))$ (Kostant's theorem, see \cite{D}). Let $V^k\subset U(su(2)_\h)$, $k\ge 0$ (we set $V^0\equiv \R$),
be an irreducible adjoint $U(su(2)_\h)$-submodule such that $b^k=(-ix-y)^k$ is a highest weight element in its complexification
$V^k_{\C}$. Thus, $V^k$ is the submodule of the spin $k$,  $\dim V^k=2k+1$.

Then we have
$$
U(su(2)_\h)= \bigoplus_{k=0}^\infty \left( Z(U(su(2)_\h))\ot V^k\right).
$$
As is known, $Z(U(su(2)_\h))$ consists of all polynomials in $\Cas$. By fixing a basis $\{v_i^k,\, 1\leq i\leq 2k+1\}$
in each submodule $V^k$ and the basis $\{\Cas^p\}$, $p\ge 0$ in $Z(U(su(2)_\h))$, we get a basis $\{\Cas^p \ot v_i^k\}$ of
the whole algebra $U(su(2)_\h)$.

A similar decomposition hold true for the algebras $U(u(2)_\h)$ and $\Ah$:
$$
U(u(2)_\h)= \bigoplus_{k=0}^\infty \left( Z(U(u(2)_\h))\ot V^k\right),\qquad
\Ah= \bigoplus_{k=0}^\infty \left( Z(\Ah)\ot V^k\right).
$$

Consider now a natural action of the rotation group $SO(3)$ on the algebra $\W(\Ah)$ (and consequently on $\Diff$). This group
acts on the "spacial generators" $x, y, z$ of the algebra $\W(\Ah)$ and on the corresponding derivatives $\pa_x$, $\pa_y$ and  $\pa_z$
in the classical  way, but it affects neither the factor $\K(\tau,\rh)$ nor the derivative $\partial_\tau$.

\begin{definition}\rm
\label{three} We call a differential operator $\DD_A\in \Diff$ to be invariant if it is the image of an $SO(3)$-invariant element
$A\in \W(\Ah)$ of the algebra $\W(\Ah)$ under the isomorphism $\W(\Ah)\rightarrow {\Diff}$.
\end{definition}

Important examples of the first order invariant operators are
\be
\tilde\partial_t \qquad {\rm  and} \qquad Q =x\, \pa_x+y\,\pa_y+z\,\pa_z.
\label{first}
\ee
Some the second order differential operators to be used below are
\be
\De_0=\dd^2,\qquad \De_1=\De=\pa_x^2+\pa_y^2+\pa_z^2,\qquad \De_2=Q\,\dd,\qquad \De_3=Q^2.
\label{second}
\ee
All operators from the lists (\ref{first}) and (\ref{second}) will be called {\em basic}. The operator $\De$ is called {\em Laplacian}.
All polynomials in the basic operators with coefficients from $\K(\tau,\rh)$ are also invariant operators.

\begin{proposition}
Let $\DD: \Ah\to\Ah$ be a basic invariant differential operator. Then it maps the component
\be
{\V}^k=\K(\tau,\rh)\ot V^k
\label{isoco}
\ee
to itself. Moreover, its action on any element $f(\tau, \rh)\, v$, $v\in V^k$ alters only the central function
$f(\tau, \rh)$. Thus, we have
$$
\DD(f(\tau, \rh)\, v) = g(\tau, \rh)\, v
$$
for some $g(\tau, \rh) \in \K(\tau, r_\hbar)$.
\end{proposition}

{\bf Proof.} First, observe  that  the {\em evaluation map}
$$
\Diff\ot \Ah\to \Ah,\,\, \DD\ot n\mapsto \DD(n)
$$
is $SO(3)$-covariant. This entails that if $\DD$ is an invariant differential operator  it commutes with the action of the group
$SO(3)$. Therefore, its complexification commutes with the action of the group $SL(2,\C)$.

Now we use the following fact, proved in \cite{GS2}.
\begin{lemma}
\label{L3}
The action of the partial derivatives $\pa_t,\dots,\pa_z$ onto the elements $b^k\in V^k_{\C}$ is subject to the classical
Leibniz rule.
\end{lemma}

Due to lemma \ref{L3} we can conclude that all basic differential operators act on the elements  $b^k$ by the classical formula.
Namely, we have
$$
\begin{array}{l}
\pa_t(b^k)=0,\qquad Q(b^k)=k\, b^k,\\
\rule{0pt}{7mm}
\displaystyle
\De(b^k)=0,\qquad \De_0(b^k)=\left( \frac{2}{h}\right)^2 b^k,\qquad \De_2(b^k)= \frac{2k}{\h}\, b^k,\qquad
\De_3(b^k)=k^2\, b^k.
\end{array}
$$

Now, the claim of the proposition follows immediately from the invariancy of the basic operators.\hfill\rule{6.5pt}{6.5pt}
\medskip

As follows from the above consideration, the restriction of any basic operator to the components ${\V}^k$ can be well defined. Note that being restricted to the component
${\V}^k$ such an operator acts as a difference operator  on the factor $f(\tau,\rh)$ in an element $f(\tau,\rh)v$, $v\in V^k$. The reader
is referred to \cite{GS2} for explicit examples and detail.

Thus, the spectral problem for a basic differential operator $\DD$
$$
\DD(f)=\la  f,\qquad f\in \Ah
$$
reduces to the spectral problem of a family of its restrictions to the components ${\V}^k$.

\begin{definition}\rm
The restriction of a basic operator $\DD$ on the component $\VV^0=K(\tau, \rh)$  is called {\it the radial part}
of the operator $\DD$ and is denoted $\DD_{rad}$.
\end{definition}

As follows from \cite{GS2}, the radial part of the Laplacian $\De$ has the form:
\begin{eqnarray}
\De_{rad}(f(\tau,\rh))=\frac{f(\tau+\hh, \rh+2\hh)+f(\tau+\hh, \rh-2\hh)-2f(\tau+\hh, \rh)}{4\hh^2} \\   \nonumber
+\frac{1}{\rh}\frac{f(\tau+\hh, \rh+2\hh)-f(\tau+\hh, \rh-2\hh)}{2\hh}
\label{radD}
\end{eqnarray}
Another important example is the radial part of the operator $Q$
\be
Q_{rad}(f(\tau, \rh))=\frac{(\rh^2-\hh^2)}{\rh}\frac{(f(\tau+\hh,\rh+\hh)-f(\tau+\hh,\rh-\hh))}{2\hh}.
\label{radQ}
\ee
This result immediately follows from (\ref{pax}) and its $y$- and $z$-analogs.

\section{Quantum radial derivative and its applications}

Let us introduce the operator $\partial_{r_\hbar}: \K(\tau,r_\hbar) \rightarrow \K(\tau,r_\hbar)$ by the following
rule
$$
\pa_{\rh}f(\tau,\rh)=\frac{f(\tau+\hbar,\rh+\hh)-f(\tau+\hbar,\rh-\hh)}{2\hh},\qquad \forall\,f(\tau,r_\hbar)\in\K(\tau,r_\hbar).
$$
This operator is well defined on the component $\K(\tau, \rh)$ and we call it {\em the quantum radial derivative}.
It is a difference operator which turns into the usual derivative in $r$ as $\hh\to 0$:
$$
\lim_{\hbar\rightarrow 0}\pa_{\rh}f(\tau,\rh) = \frac{\partial f(\tau, r)}{\partial r}.
$$
With this notation we rewrite formula (\ref{pax}) in the form:
\be
\pa_x(f(\tau,\rh))=\frac{x}{\rh}\,\parh(f(\tau, \rh)). \label{xcase} \ee
This formula is similar to the classical one but with another meaning of  the derivative.

Formula (\ref{radQ}) can also be rewritten in terms of the new derivative:
$$
Q_{rad}(f(\tau, \rh))=\frac{(\rh^2-\hh^2)}{\rh}\,\parh(f(\tau, \rh)).
$$
Since the function $f(\tau, r_\hbar)$ is an arbitrary element of $\K(\tau, r_\hbar)$ we have an operator
equality on the algebra $\K(\tau, r_\hbar)$
\be
\parh=\frac{\rh}{\rh^2-\hh^2}\,Q.
\label{parh}
\ee
We use this relation as the definition of an extension of the quantum radial derivative $\partial_{r_\hbar}$ on the whole
algebra $\Ah$. Note that this definition differs from that from \cite{GS2} by the factor at the operator $Q$.

As the next step we express the operator $\De_{rad}$ in terms of  the quantum radial derivative.
A straightforward calculation on the base of relation (\ref{radD}) leads to the result:
\be
\De_{rad}(f(\tau,\rh))=\frac{1}{\rh}\,\parh^2(\rh\,f(\tau, \rh)).
\label{rrLap}
\ee

\begin{remark} \rm
Formula (\ref{rrLap}) is analogous to the classical one
\be
\De_{rad}(f(t,r))=\frac{1}{r}\partial^2_r(rf(t,r)) = \pa_r^2(f(t,r))+\frac{2}{r}\pa_r(f(t,r))
\label{rLap}
\ee
and turns into this as $\hh\to 0$.
Actually, formula (\ref{rrLap}) has motivated our definition of the quantum radial derivative.

In contrast with the classical case, the Laplace operator $\De$ affects the variable $\tau$,
though it does not explicitly contain the derivative $\pat$. This effect disappears as $\hh\to 0$.
\end{remark}

Our next aim is to solve the equation
\be
\De(f(\rh))= \De_{rad}(f(\rh))=0
\label{solv}
\ee
for a $\tau$-independent function $f(\rh)\in Z(\Ah)$. In the classical case ($\hh=0$) the vector space of solutions of the equation $\De_{rad}(f(r))=0$ in the region $r>0$
consists of all linear combinations of two functions $f_1(r)=1$ and $f_2(r)=\frac{1}{r}$.

In the quantum case as follows from expression (\ref{rrLap}) the functions $f_1(\rh)=1$ and $f_2(\rh)=\rh^{-1}$
are solutions of the equation (\ref{solv}). Moreover, all linear combinations $C_1(\rh)\, f_1+C_2(\rh)\, f_2$ are also
solutions provided that the coefficients $C_i(\rh)$, $i=1,2$, are $2\hh$-periodical functions in $\rh$.

\begin{remark}
\label{onze}
\rm
In the classical case the function $\frac{1}{r}=\frac{1}{\sqrt{x^2+y^2+z^2}}$ satisfies the equation
$$
\De(\frac{1}{r})=-4{\pi}{\de}(\mbox{\bf{x}})
$$
on the whole space $\R^3$ . Here, $\mbox{\bf{x}}=(x,y,z)$ and $\de(\mbox{\bf{x}}-\mbox{\bf{a}})$ stands for
the $\de$-function located at the point $\mbox{\bf{a}}= (a_1,a_2, a_3)\in \R^3$. In section 6 we shall discuss
a quantum version of this equation.
\end{remark}

Note that  the operator (\ref{rrLap})  is well defined on a lattice
\be
(\rh, \tau)\in \{a+2k\hh, b+l\hh,\,\, k,l\in \Bbb{Z}\},
\label{lat}
\ee
where $a$ and $b$ are some fixed real numbers.
Thus, we can consider the restriction of this operator to the lattice  (\ref{lat}). A similar restriction is possible for any operator
depending only on the quantities $\tau$ and $\rh$ and the  derivatives $\pat$ and $\parh$. In this sense we speak about the
{\em discretization} of a dynamical model described by such an operator.

This discretization reduces the freedom for coefficients in the space of solutions of the equation  $\De(f(\tau, \rh))=0$. Thus,
the aforementioned coefficients  $C_i(\rh)$ entering  the general solution become constant  on each of the two
 sublattices
$$
\{a+2k\hh,\,\,b+l\hh,\, k,l\in \Bbb{Z}\}\quad {\rm and}\quad \{a+\hh+2k\hh,\,\,b+l\hh,\, k,l\in \Bbb{Z}\}.
$$

Let us assume the algebra $U(gl(2))$ to be represented in a Verma module $M_\La$ of the highest weight $\La=(\La_1,\La_2)$, $\La_1\geq \La_2$
(for the definition we refer the reader to the book {\cite{M}).  The module $M_\La$ has a finite dimensional $U(gl(2))$-submodule iff $\La_1-\La_2=n$ is
a non-negative integer. On passing to the algebra $U(u(2))$ we get (after rescaling of the generators $x\rightarrow x/h$ etc.) the corresponding
representation $\rho_\La: U(u(2)_h)\to \End(M_\La)$, which can also be factorized up to a finite dimensional representation provided $\La_1-\La_2=n$.
(Note that $\frac{n}{2}$ is the spin of this finite dimensional representation.)

It is known that the image of any central element $u\in Z(U(u(2)_h))$ under the representation $\rho_\La$ is a scalar operator (it will be identified with
its eigenvalue). In particular, we have\footnote{In this connection we would like to mention the papers \cite{Gr, NT} where a CH identity (and its image
in Verma and finite dimensional modules) was considered. As follows from the cited papers (also, see \cite{M}) the roots of this CH identity, after a
proper ordering, are connected with the weight $\La$ by a relation $\mu_1=\La_1+1$, $\mu_2=\La_2$ (for $h=1$).}
$$
\rho_\La(\Cas)=(\La_1-\La_2)(\La_1-\La_2+2)\hh.
$$
 This enables us to compute the corresponding value of the quantum radius $\rh$. Namely, we have
$$
\rho_\La(\rh)=\sqrt{(\La_1-\La_2+1)^2\hh^2}=(\La_1-\La_2+1)\hh.
$$
The sign  here is motivated by our wish to have a positive quantum radius (at $\hbar >0$) on finite dimensional modules.
Thus, on the finite dimensional module of the spin $\frac{n}{2}$ the quantum radius takes the values
$$
\rh=(n+1)\hh, \quad n\in{\Bbb Z}_+.
$$
So, the operator $\Delta_{rad}$ admits a discretization well defined on a family of the Verma modules, some of them can be factorized to finite dimensional ones.
This is also true for the isotypic components ${\cal V}^k$ (\ref{isoco}) with $k\not=0$.

We complete this section with the following observation. Let
\be
\al: \Sym(u(2))\to U(u(2)_\h)
\label{qmap}
\ee
be an $SU(2)$-covariant quantizing map (see \cite{GS1} for detail). Then any operator
$$
\D:\Sym(u(2))\to \Sym(u(2))
$$
can be pulled forward to
\be
\D_\al: U(u(2)_\h)\to U(u(2)_\h),\qquad \D_\al=\al \, \D\, \al^{-1}.
\label{D-al}
\ee
 However, the image of a usual partial derivative in the algebra $\Sym(u(2))$ with respect to the map (\ref{D-al}) does not
coincide with the quantum partial derivatives considered in the present paper.  Nevertheless, the images of the infinitesimal rotations in
$\Sym(u(2))$ remain to be infinitesimal rotations in $U(u(2)_\h)$.

As for the boosts, their analogs are not naturally defined on the algebra $U(u(2)_\h)$. Let us introduce the infinitesimal  "radial boost"
in the classical case $(\hbar = 0)$
\be
t\pa_r+r\pa_t=t\frac{Q}{r}+r\pa_t.
\label{boo}
\ee
The corresponding one parameter group $\exp(\nu (t\pa_r+r\pa_t))$ preserves the Lorentz interval
 $s =\sqrt{t^2-(x^2+y^2+z^2)}=\sqrt{t^2-r^2}$ since
$$
(t\pa_r+r\pa_t)(s) = 0.
$$

We define the quantum analog of the "infinitesimal radial boost" as
$$
\tau\parh+\rh\pa_\tau
$$
and the "quantum Lorentz interval" squared as $\tau^2-\rh^2$. However,  now the action
$$
(\tau\parh+\rh\pa_\tau)(\tau^2-\rh^2)=\tau(-2\rh)+\rh(2\tau-2\hh)=-2\hh\rh
$$
is proportional to $\hh$ and disappears only as $\hh\to 0$.

\section{Klein-Gordon and  Schr\"odinger equations on $U(u(2)_\hbar)$}

In this section we consider the NC (quantum)  versions of two dynamical models playing an important role in physics.
Namely, we are dealing with the Klein-Gordon equation (as an example of a relativistic model) and the Schr\"odinger equation
with the hydrogen atom potential (as an example of a non-relativistic model). For the former equation we exhibit an analog of
the plane wave solution under an assumption on momenta parameterizing the solution. With this assumption we get a deformed
version of the relativistic dispersion relation. For the latter equation we are looking for a value of the ground state energy.

\medskip

\leftline{\bf The Klein-Gordon equation}

The Klein-Gordon equation describes evolution of a free massive scalar field. It is the second order partial differential
equation of the form (we assume the speed of light and the Planck constant to be  equal to the unity)
\be
(\partial_t^2-\partial_x^2-\partial_y^2-\partial_z^2 +m^2)(\phi(t,x,y,z)) = 0.
\label{eq:Kl-G}
\ee
Here $m$ is a real constant (the mass of the field). This equation possesses the plane wave solutions
\be
\phi(t,x,y,z) = \exp i(Et-p_xx-p_yy-p_zz)
\label{plw} \ee
where the energy $E$ and the momenta $(p_x,p_y,p_z)$ are subject to the relativistic dispersion relation
$$
E^2 = p_x^2 + p_y^2 + p_z^2 + m^2.
$$

Let us pass to a NC version of the Klein-Gordon equation. Since the generators $x,y, z$ do not explicitly enter
equation (\ref{eq:Kl-G}), this passage is straightforward --- we have to assign a new meaning to the derivatives expressed
in the relations\footnote{Here we use the new system of generators and parameters $\tau = -it$ and  $\hbar = h/2i$.} (\ref{leib-r}).

Note that in the NC case we meet an additional difficulty. The symbol $\phi(\tau,x,y,z)$ is not well defined now, since we have to
explicitly fix some ordering of the NC variables $x$, $y$ and $z$. To avoid this technical difficulty, we consider a solution depending
on the central element $\tau $ and the only spacial variable, say $x$. Otherwise stated, we assume the momentum to be of the form
$p=(p_x,0,0)$.

\begin{remark}
{\rm This simplification enables us to compute the action of the Klein-Gordon operator on a plane wave (\ref{plw}). Note that
since the element $t, x,y, z$ belong to the algebra $U(u(2)_h)$, this plane wave can be treated as  an element of the group
$U(2)$ (at least for the small enough coefficients $E$ and $p$) presented in the coordinates of the first kind. The computation
becomes more easy if we present elements of this group in the coordinates of the second kind. This means that the function
$\phi$ has the form $ \exp(iEt)\,\exp(- i p_xx)\, \exp(- i p_yy)\, \exp(- i p_zz)$. An example of computation with such a function
is given in \cite{GPS2}.}
\end{remark}

First of all, we note, that similarly to the classical case
$$
(\partial_y^2+\partial_z^2)(\phi(\tau,x)) = 0.
$$
This immediately follows from formula (\ref{ap1}) of Appendix. So, our equation becomes
\be
(\partial^2_\tau -\partial^2_x +m^2)(\phi(\tau,x)) = 0.
\label{eq:KG-1}
\ee

 We look for  solutions of this equation in the form
$$
\phi(\tau,x) = f_E(\tau)g_p(x)
$$
where $E$ and $p=p_x$ are numeric parameters of the solution (analogs of the energy and momentum of the commutative
case). Our goal is to find the explicit form of the functions $f_E(\tau)$ and $g_p(x)$ and the dispersion relation, connecting
$E$ and $p$.

Let us demand the function $f_E(\tau)$ to be an eigenfunction of the derivative $\partial_\tau$ with the eigenvalue $E$.
The action of the derivative $\partial_\tau$ on a function in $\tau$ can be extracted from (\ref{fo}). So, we have
$$
\partial_\tau (f_E(\tau)) = \frac{f_E(\tau +\hbar) - f_E(\tau)}{\hbar} = Ef_E(\tau).
$$
The solution of this difference equation reads:
\be
f_E(\tau) = \xi(E)^{\tau/\hbar},\quad {\rm where}\quad \xi(E)=1+\hbar E.
\label{sol}
\ee
Note that the function $f_E(\tau)$ can be presented in the form $f_E(\tau)=\exp(\al(E) \tau)$, where
$\al(E)=\frac{\ln (1+\hh E)}{\hh}$. In the classical  limit $\hbar\rightarrow 0$ the function $\al(E)$ turns into $E$.

In a similar way, the function $g_p(x)$ is taken to be an eigenfunction of the operator $\partial_x$ with the eigenvalue $p$:
$$
\partial_x (g_p(x)) = \frac{g_p(x+\hbar) - g_p(x-\hbar)}{2\hbar} = p\,g_p(x).
$$
We take the solution of this equation in the form
$$
g_p(x) = \eta(p)^{x/\hbar}, \quad {\rm where}\quad \eta(p)=\hbar p+\sqrt{1+(\hbar p)^2}.
$$
In the classical limit we get $\lim_{\hbar\rightarrow 0} g_p(x)= e^{xp}$.

Let us now turn to solving equation (\ref{eq:KG-1}) with the above substitution $\phi(\tau,x) = f_E(\tau)g_p(x)$.
First of all, permuting $\partial_\tau^2$ and $f_E(\tau)$ we have:
$$
\partial_\tau^2 f_E(\tau) = E^2f_E(\tau) +2Ef_E(\tau +\hbar)\partial_\tau+f_E(\tau+2\hbar)\partial_\tau^2
=f_E(\tau)(E^2 +2E\xi \partial_\tau +\xi^2\partial_\tau^2).
$$
Also (see Appendix for detail)
$$
\partial_\tau (g_p(x)) = g_p(x)\, \frac{(\eta-1)^2}{2\hbar\,\eta},\qquad \partial_\tau^2 (g_p(x))  = g_p(x)\,\frac{(\eta-1)^4}{(2\hbar\,\eta)^2}.
$$
In the above formulae we omit the arguments $E$ and $p$ in the notations $\xi(E)$ and $\eta(p)$ respectively. All these formulae follow
from the permutation relations of the partial derivatives in $\tau$,   presented in Appendix, and  the particular form of the
functions $f_E(\tau)$ and $g_p(x)$. This allows us to calculate the term $\partial^2_\tau (\phi(\tau,x))$ in (\ref{eq:KG-1}).

In order to find the value of the second term $\partial^2_x(\phi(\tau,x))$, we use  the relation
$$
\partial_x^2 (f_E(\tau)g_p(x)) = f_E(\tau+2\hbar)\partial_x^2(g_p(x)) = p^2\xi^2f_E(\tau)g_p(x).
$$
Finally, we come to the  equation:
$$
f_E(\tau)g_p(x)\left[ \Bigl(\xi\,\frac{(\eta-1)^2}{2\hbar\,\eta} +E\Bigr)^2 -\xi^2p^2+m^2\right] =0,
$$
which, in turn, leads to the dispersion relation, connecting the parameters $E$ and $p$
$$
\Bigl(\xi(E)\,\frac{(\eta(p)-1)^2}{2\hbar\,\eta(p)} +E\Bigr)^2 -\xi(E)^2p^2+m^2 = 0.
$$

On taking into account the explicit form of $\xi(E)$ and $\eta(p)$, we simplify the dispersion relation to the form
$$
E^2-2p^2\,\frac{1+\hbar E}{1+\sqrt{1+(\hbar p)^2}}+m^2 = 0.
$$
On solving the equation in $E$ we find
$$
E_\pm = \pm\sqrt{p^2-m^2} +\frac{\hbar p^2}{1+\sqrt{1+(\hbar p)^2}}.
$$
\medskip

This is the final form of the "quantum dispersion relation" under the above assumption $p=p_x$, $p_y=p_z=0$.

\medskip

\leftline{\bf The Schr\"odinger type equation}
Now, let us consider a differential equation of the Schr\"odinger type
\be
(a\partial_\tau+b\Delta+\frac{q}{\rh})\,(\psi) = 0
\label{eq-schr}
\ee
where $r_\hbar$ is a quantum radius and $a$, $b$, $q$ are some numeric constants. It is convenient to look for  $\psi$ in
the form
$$
\psi=\psi^{(k)}(\tau,\rh)\, v = f_E(\tau)\phi(\rh)\, v,
$$
where $v\in V^k$ and $\psi^{(k)}(\tau,\rh)\in Z(\Ah)$ and $f_E(\tau)$ is defined by formula (\ref{sol}).

Below we are dealing with the particular case  $k=0$ of this problem. Using the properties of $f_E(\tau)$ we get the equation on $\phi(\rh)$:
\begin{eqnarray*}
aE\phi(\rh)&+& a\xi(E)\,\frac{\phi(\rh+\hbar)(\rh+\hbar)+\phi(\rh-\hbar)(\rh-\hbar)-2\rh \phi(\rh)}{2\hbar\,\rh}\\
&+& b\xi(E)^2\left[\frac{\phi(\rh+2\hbar)-\phi(\rh-2\hbar)}{2\hbar\,\rh}+\frac{\phi(\rh+2\hbar)
+\phi(\rh-2\hbar)-2\phi(\rh)}{4\hbar^2}\right]\\
&+& \frac{q}{\rh}\,\phi(\rh)=0,\\
\end{eqnarray*}
where $\xi(E)$ is defined in (\ref{sol}).

Now, we are interested in question: is there a ground state of the form $\phi(\rh) = e^{- \sigma \rh}$ (similarly to the commutative
case)? Substituting this exponent into the equation, we get
\begin{eqnarray*}
aE\rh&+& \frac{a\xi(E)}{2\hbar}\,\left(\rh(e^{\sigma\hbar}+e^{-\sigma\hbar}-2)-\hbar(e^{\sigma\hbar}-e^{-\sigma\hbar}) \right)\\
 &+& \frac{b\xi(E)^2}{4\hbar^2}\left[ \rh(e^{2\sigma\hbar}+e^{-2\sigma\hbar}-2)-2\hbar(e^{2\sigma\hbar}-e^{-2\sigma\hbar})  \right]+q=0.\\
\end{eqnarray*}
Collecting separately the terms with  $\rh$ and those without it, we get a system of equations to determine $\sigma$ and $E$:
$$
\left\{
\begin{array}{l}
\displaystyle aE+2\,\frac{a\xi(E)}{\hbar}\,\sinh^2(\sigma\hbar/2)+\frac{b\xi^2(E)}{\hbar^2}\,\sinh^2(\sigma\hbar) = 0\\
\rule{0pt}{8mm}
\displaystyle a\xi(E)\sinh(\sigma\hbar)+\frac{b\xi^2(E)}{\hbar}\,\sinh(2\sigma\hbar)-q = 0
\end{array}\right.
$$
or
\be
\left\{
\begin{array}{l}
\displaystyle \xi(E) = 2\,\frac{\cosh(\sigma\hbar)-\sinh(\sigma\hbar)q/(2a)}{1+\cosh^2(\sigma\hbar)}\\
\rule{0pt}{8mm}
\displaystyle \xi^2(E) = \frac{a\hbar}{b}\,\frac{1-\xi(E)\cosh(\sigma\hbar))}{\sinh^2(\sigma\hbar)}
\end{array}.
\right.
\label{gr-energy}
\ee
This system reduces to the equation on $y=\mathrm{tanh}(\sigma\hbar)$:
\be
(\omega/2-2\rho^2) y^3+\rho(4-\omega)y^2 - (\omega +2)y+2\rho\omega = 0,
\label{eq-tha}
\ee
where
$$
\rho=\frac{q}{2a},\quad \omega=\frac{a\hbar}{b}.
$$
In the commutative limit $\hbar\rightarrow 0$ the above equation gives in the leading order in $\hbar$
\be
\sigma_0 = \frac{q}{2b},
\label{lim-a}
\ee
while the leading order of the ground state energy can be extracted from the first line of the system
(\ref{gr-energy})
\be
E_0 = -\,\frac{q}{2a}\,\sigma_0 = -\,\frac{q^2}{4ab}.
\label{lim-E}
\ee
These limit values coincide with the corresponding solution of the differential equation
$$
\left(aE +b\left(\frac{2}{r}\,\partial_r+\partial^2_r\right) +\frac{q}{r}\right)\,(\phi(r)) = 0,
$$
which is the commutative limit of (\ref{eq-schr}) with zero value of the orbital momenta.
Bellow, we compute the first order corrections in $\hbar$ to the values $E_0$ and $\sigma _0$.
We do not know any "physical meaning" of other solutions of the equation (\ref{eq-tha}).

The first order corrections in $\hbar$ to the classical results (\ref{lim-a}) and (\ref{lim-E}) can be found from the
system (\ref{gr-energy}). We present the result of calculations by passing to the physical meaning of the parameters $a,\,b,\,q$ entering the equation (\ref{eq-schr}). More precisely, we set
$a=h_{pl} c$, $b=h^2_{pl}/2m$, $q= e^2$, where $h_{pl}$ is the Planck constant, $c$ is the speed of light, $m$ is the
electron mass, $e^2$ is the electron charge squared. Besides, we must renormalize the energy parameter
$E\rightarrow E/(h_{pl}c )$. It is convenient to introduce the dimensionless fine structure constant $\alpha$ and
the electron De Broglie wave length $\lambda_B$ by the relations
$$
\alpha = \frac{e^2}{h_{pl}c},\qquad \lambda_B = \frac{h_{pl}}{mc}.
$$
With this constants our answer has the following simple form:
\begin{eqnarray}
&&E = -mc^2\,\frac{\alpha^2}{2}\left(1-\frac{\hbar}{4\lambda_B}\,\left(1-\frac{\alpha^2}{2}\right) \right)+O(\hbar^2)\label{E-first}\\
&&\rule{0pt}{7mm}
\sigma = \frac{\alpha}{2\lambda_B}\,\left(1-\frac{\hbar}{\lambda_B}\,(1-\alpha^2)\right)+O(\hbar^2).\label{sig-first}
\end{eqnarray}
The above relations show that noncommutativity of the space-time leads to increasing of the ground state energy as well as
to increasing of an average size of the hydrogen atom which is roughly proportional to $\sigma^{-1}$ (we assume $\hbar$ to be positive).

As for the eigenfunction $f_E(\tau)e^{-\sigma r_\hbar}$, the first order correction to it can be easily obtained with the use
of (\ref{E-first}) and (\ref{sig-first}) and is left for the reader.

\section{Concluding remarks, discussion}

1. First, let us go back to the remark \ref{onze}. The relation exhibited there means that for any sufficiently
smooth compactly supported  function $f(x,y,z)$ defined on $\R^3$ we have
\be
\int_{\R^3} \De(\frac{1}{r})\,f(x,y,z) \, dxdydz=-4\pi\, f(0),
\label{int}
\ee
where the above integral is treated in the sense of the principal value.

A proof of (\ref{int}) is usually based on the  Green identities. We slightly modify  this proof in order to adapt
it to the quantum case.

Since we consider the function $\frac{1}{r}$ as a distribution, the integral in the left hand side of (\ref{int}) is
by definition
\be
\lim_{\varepsilon\to 0, R\to \infty}\int_{V(R,\varepsilon)} \frac{1}{r}\,\De(f(x,y,z)) \, dxdydz,
\label{intt}
\ee
where $V(R,\varepsilon)=B(0,R)\setminus B(0,\varepsilon)$ and $B(0,a)$ is the ball with the center at 0 and
the radius $a\geq 0$.

Note that it suffices to assume the function $f$ to be dependent on the radius $r$ only. Otherwise, it can
be replaced by
$$
g(r)=\frac{1}{4\pi}\int f(x,y,z) d\Om
$$
where $d\Om$ is the surface 2-form of the unit sphere.  Note that $g(0)=f(0,0,0)$.

Then, in virtue of (\ref{rLap}), the integral (\ref{intt}) can be written as
$$
4\pi\,\int_{\varepsilon}^{R}\frac{1}{r}\left(\frac{1}{r}\pa_r^2(r f(r))\right)\, r^2 dr=4\pi\,(f(R)+Rf'(R)-f(\varepsilon)-
\varepsilon f'(\varepsilon)).
$$
On passing to the limits $R\to \infty,\, \varepsilon \to 0$ and assuming the function $f$ to be continuously
differentiable at $(0,0,0)$ we get the claim.

Heuristically, we can  reproduce the computation of (\ref{intt}) in the quantum case (using of course (\ref{rrLap})
instead of (\ref{rLap})) and get the same result under the following assumptions (actually, definitions)
$$
dxdydz=d\Om\, \rh d\rh\, \qquad {\rm and}\qquad \int_a^b \parh(f(\rh))\, d\rh=f(b)-f(a).
$$

Thus,  the notion of the integral (including that "in spherical coordinates") can be adapted to the derivative $\parh$
in the classical manner. We plan to consider a multidimensional analog of this formula related to the algebra
$U(gl(m)_h)$ in the subsequent publications.

2. Note that we have also considered another basic physical model, namely, 3-dimensional harmonic oscillator. Compared with
the hydrogen atom model, the harmonic oscillator is more difficult to deal with. The main difficulty here is in the asymptotic  of
the wave function of the commutative case $\phi(r)\sim e^{-\sigma r^2}$. The direct substitution of the non-commutative
analog $e^{-\sigma r_\hbar^2}$ into the corresponding equation leads to an equation which we cannot solve.
Apparently, one should find here some another function in $r_\hbar$ (with the commutative limit $e^{-\sigma r^2}$) which
would be more adequate to the problem in the NC case. Another possible way out is a modification of the harmonic oscillator
potential (with recovering its classical form at the commutative limit $\hbar \rightarrow 0$) in such a way that the exponent
$e^{-\sigma r_\hbar^2}$ would be an eigenfunction of the modified Hamiltonian.

3. The renormalization of the time variable $\tau$ leads to a change of the step of the corresponding
lattice. We do not know what normalization of the time variable and the corresponding step could be "physically motivated".

4. Now, let us compare once more our approach to quantization of dynamical models with that arising from deformation quantization.
As we noticed in Introduction, the models introduced in the frameworks of the latter approach
make use of the classical derivatives. Only, the usual commutative product of functions on a given Poisson variety
is  replaced by the corresponding $\star$-product. Note that sometimes such a $\star$-product is
introduced via a map  (\ref{qmap}). For instance, this map can be constructed with the help
 of isotypic components as done in \cite{GS1}. In this case it can be naturally extended to the completion similar to that $\Ah$.

Then the operator describing a given
dynamical model can be pulled forward to the quantum algebra by means of this map (see formula (\ref{D-al})).
However, this way of
quantizing dynamical models is equivalent to the previous one and does not give rise to partial derivatives which could be difference operators. Consider the operator $\parh$ defined in this way on the center of the algebra $\Ah$
(see \cite{GS1}).  Since this operator is the image of the derivative $\pa_r$
under the map $f(r)\to f(\rh)$, it   is also a
differential operator. By contrast, the operator $\parh$ constructed in the current paper is a difference operator.

\section{Appendix}
Here we collect some technical results on quantum differential calculus.
First, observe that in the new system of quantities  $\tau=-i t$,
$\hbar=\frac{h}{2i}$, $\rh=\frac{\mu}{2i}$ the generating system of the algebra $U(u(2)_h)$ is
$$
[x, \, y]=2i\hbar\, z,\quad [y, \, z]=2i\hbar\, x,\quad[z, \, x]=2i\hbar\, y,\quad [\tau, \, x]=[\tau, \, y]=[\tau, \, z]=0.
$$
whereas the permutation relations between the generators and derivatives are
$$
\begin{array}{l@{\qquad}l@{\qquad}l}
[\partial_\tau,\tau] = 1+\hbar\,\partial_\tau, & [\partial_x,\tau] = \hbar\,\partial_x,&[\partial_x,y] = i\hbar\,\partial_z,\\
\rule{0pt}{7mm}
[\partial_\tau ,\,x]= \hbar\,\partial_x, &[\partial_x,x] = 1+ \hbar\,\partial_\tau, & [\partial_x,z] = -i \hbar\,\partial_y.\\
\end{array}
$$
The other relations can be obtained from these by the cyclic permutations $x\rightarrow y\rightarrow z$.

It is convenient to work with the "shifted derivative" $\tilde \partial_\tau = \partial_\tau+1/\hbar$. Thus, we get
$$
\tilde\partial_\tau f(\tau) = f(\tau +\hbar)\tilde\partial_\tau\qquad\mathrm{or}\qquad
\partial_\tau f(\tau) = \frac{f(\tau+\hbar)-f(\tau)}{\hbar} +f(\tau+\hbar)\partial_\tau.
$$
Next, the permutations of the derivatives $\tilde\partial_\tau$ and $\partial_x$ with $x$ can be cast in the matrix form
$$
\left(\!\!\begin{array}{c}
\tilde\partial_\tau\\
\partial_x
\end{array}\!\!\right)x =
\left(\!\!\begin{array}{cc}
x&\hbar\\
\hbar& x
\end{array}\!\!\right)\left(\!\!\begin{array}{c}
\tilde\partial_\tau\\
\partial_x
\end{array}\!\!\right) = \left[\frac{x+\hbar}{2}
\left(\!\!\begin{array}{cc}
1&1\\
1& 1
\end{array}\!\!\right)+\frac{x-\hbar}{2}
\left(\!\!\begin{array}{cc}
1&-1\\
-1&1
\end{array}\!\!\right)\right]\left(\!\!\begin{array}{c}
\tilde\partial_\tau\\
\partial_x
\end{array}\!\!\right).
$$

This relation enables us to get
\begin{eqnarray*}
\tilde\partial_\tau f(x) &=& f(x+\hbar)\,\frac{\tilde\partial_\tau+\partial_x}{2} +f(x-\hbar)\,\frac{\tilde\partial_\tau-\partial_x}{2},\\
\partial_x f(x) &=& f(x+\hbar)\,\frac{\tilde\partial_\tau+\partial_x}{2} +f(x-\hbar)\,\frac{\partial_x-\tilde\partial_\tau}{2},
\end{eqnarray*}
or
\begin{eqnarray*}
\partial_\tau f(x) &=& \frac{f(x+\hbar)+f(x-\hbar)-2f(x)}{2\hbar}+f(x+\hbar)\,\frac{\partial_\tau+\partial_x}{2} +f(x-\hbar)\,\frac{\partial_\tau-\partial_x}{2},\\
\partial_x f(x) &=& \frac{f(x+\hbar)-f(x-\hbar)}{2\hbar}+f(x+\hbar)\,\frac{\partial_\tau+\partial_x}{2} +f(x-\hbar)\,\frac{\partial_x-\partial_\tau}{2}.
\end{eqnarray*}

In an analogous way we get
\begin{eqnarray*}
\partial_y f(x) &=& f(x+\hbar)\,\frac{\partial_y-i\partial_z}{2} +f(x-\hbar)\,\frac{\partial_y+i\partial_z}{2},\\
\partial_z f(x) &=& f(x+\hbar)\,\frac{\partial_z+i\partial_y}{2} +f(x-\hbar)\,\frac{\partial_z-i\partial_y}{2},
\end{eqnarray*}
or
\begin{eqnarray*}
(\partial_y +i\partial_z)f(x) &=& f(x-\hbar)(\partial_y +i\partial_z),\\
(\partial_y -i\partial_z)f(x) &=& f(x+\hbar)(\partial_y -i\partial_z).
\end{eqnarray*}

Consequently, we arrive to the following formulae
\be
(\partial_y^2 +\partial_z^2)f(x) = f(x)(\partial_y^2 +\partial_z^2)\quad\Rightarrow\quad \Delta(f(x)) = \partial_x^2(f(x)).
\label{ap1} \ee
Besides,
$$
\Delta f(\tau) = f(\tau+2\hbar)\Delta.
$$
Now, exhibit some other formulae containing  the second derivatives:
$$
\partial_\tau^2 f(\tau) =\frac{f(\tau+2\hbar)+f(\tau) -2f(\tau+\hbar)}{\hbar^2} +
2 \frac{f(\tau+2\hbar)-f(\tau+\hbar)}{\hbar}\partial_\tau +f(\tau+2\hbar)\partial_\tau^2,
$$
$$
\partial_\tau^2(f(x)) = \frac{f(x+2\hbar)-4f(x+\hbar)+6f(x)-4f(x-\hbar)+f(x-2\hbar)}{4\hbar^2},
$$
$$
\partial_x^2(f(x)) = \frac{f(x+2\hbar)+f(x-2\hbar)-2f(x)}{4\hbar^2}.
$$

Also, we need  permutation relations of the derivatives and functions in $\rh$:
$$
\tilde\partial_\tau f(\rh) = \frac{f(\rh+\hbar)}{2\rh}\,((\rh+\hbar)\tilde\partial_\tau +Q)+
\frac{f(\rh-\hbar)}{2\rh}\,((\rh-\hbar)\tilde\partial_\tau -Q),
$$
where $Q = x\partial_x+y\partial_y+z\partial_z$, and
$$
\partial_x f(\rh) = \frac{f(\rh+\hbar)}{2\rh}\,((\rh+\hbar)\partial_x+x\tilde\partial_\tau+ iX)+\frac{f(\rh-\hbar)}{2\rh}\,((\rh-\hbar)\partial_x-x\tilde\partial_\tau -iX),
$$
where $X = y\partial_z-z\partial_y$. Besides, we have
\begin{eqnarray*}
iX f(\rh) &=& \frac{f(\rh+\hbar)}{2\rh}\,((\rh^2-\hbar^2)\partial_x-2\hbar x\tilde\partial_\tau-xQ+ (\rh-\hbar)iX)\\
 & &+ \frac{f(\rh-\hbar)}{2r}\,(-(\rh^2-\hbar^2)\partial_x+2\hbar x\tilde\partial_\tau+xQ+ (\rh+\hbar)iX),
\end{eqnarray*}
\begin{eqnarray*}
Q f(\rh) &=& \frac{f(\rh+\hbar)}{2\rh}\,((\rh^2-\hbar^2)\tilde\partial_\tau+(\rh-\hbar)Q)\\
 & &+ \frac{f(\rh-\hbar)}{2\rh}\,(-(\rh^2-\hbar^2)\tilde\partial_\tau+(\rh+\hbar)Q),
\end{eqnarray*}
$$
\Delta(f(\rh)) = \frac{f(\rh+2\hbar)-f(\rh-2\hbar)}{2\hbar\,\rh}+\frac{f(\rh+2\hbar)+f(\rh-2\hbar)-2f(\rh)}{4\hbar^2}.
$$

\end{document}